\documentclass[12pt]{amsart}

\usepackage{amsfonts,amsmath,amsthm,amssymb,amsrefs}
\usepackage[mathscr]{eucal}
\usepackage{hyperref,xcolor}

\newtheorem{thm}{Theorem}[section]

\newtheorem{lm}[thm]{Lemma}
\newtheorem{prp}[thm]{Proposition}
\newtheorem{fact}[thm]{Fact}
\theoremstyle{definition}
\newtheorem{df}[thm]{Definition}

\newtheorem{prb}{Problem}

\newcommand{\set}[2]{\{#1\,:\,\text{#2}\}}     % { x : ... }

\newcommand{\m}[1]{{\mathbf{\uppercase{#1}}}}  % For bold uppercase algebra names

\DeclareMathOperator{\Con}{Con}
\newcommand{\calA}{\mathscr{A}}
\newcommand{\calF}{\mathscr{F}}
\newcommand{\scrK}{\mathscr{K}}
\newcommand{\calV}{\mathscr{V}}

\newcommand{\Span}{\operatorname{span}}
\newcommand{\ra}{\rightarrow}
\newcommand{\Aut}{\operatorname{Aut}}

\newcommand{\Eil}{Eilenberg-Sch\"{u}tzenberger Problem}
\newcommand{\Ho}{\mathsf{H}}
\newcommand{\Su}{\mathsf{S}}
\newcommand{\Pd}{\mathsf{P}}

\newcommand{\Alg}[3]{\mathbf{Alg}(#1,#2;#3)}
\newcommand{\Auto}[3]{\mathbf{Auto}(#1,#2;#3)}
\newcommand{\Algzero}[3]{\mathbf{Alg}^0(#1,#2;#3)}

\newcommand{\FB}{\textup{FB}}
\newcommand{\NFB}{\textup{NFB}}
\newcommand{\INFB}{\textup{INFB}}
\newcommand{\INFBfin}{\textup{INFB}${}_\mathit{fin}$}
\newcommand{\dom}{\operatorname{dom}}

\newcommand{\Vtau}{\Omega_\tau}
\begin{document}
\bibliographystyle{siam}

\title[Algebras from finite group actions]{Algebras from  finite group actions and a question of Eilenberg and 
Sch\"{u}tzenberger}
\date{Dec 1, 2022}
\thanks{The first author gratefully acknowledges the support of the Punjab Higher Education Committee (PHEC) of Pakistan.  The second author gratefully acknowledges the support of the Natural Sciences and Engineering Research Council (NSERC) of Canada.}


\author{Salma Shaheen}
\address{Pure Mathematics Department, University of Waterloo, Waterloo,
ON N2L 3G1 Canada\\ \and Mathematics Department, Government College for Women,
Dhoke Ellahi Buksh, Rawalpindi, Pakistan}
\email{salmashaheen88@gmail.com}
%\urladdr{XXX}

\author{Ross Willard}
\address{Pure Mathematics Department, University of Waterloo, Waterloo,
ON N2L 3G1 Canada}
\email{ross.willard@uwaterloo.ca}
\urladdr{www.math.uwaterloo.ca/$\sim$rdwillar}

\begin{abstract}
In 1976 S. Eilenberg and M.-P. Sch\"{u}tzenberger posed the following
diabolical
question: if $\m a$ is a finite algebraic structure,
$\Sigma$ is the set of all identities true in $\m a$,
and there exists a finite subset $F$ of $\Sigma$ such that $F$ and $\Sigma$
have exactly the same \emph{finite} models, must there also exist
a finite subset $F'$ of $\Sigma$ such that $F'$ and $\Sigma$ have exactly
the same \emph{finite and infinite} models?  (That is, must the identities of
$\m a$ be ``finitely based"?)
% and the identities
%true in the pseudovariety
%generated by $\m a$ are finitely based [bad wording: change], does it follow that the identitites
%true in the full variety generated by $\m a$ are also finitely based?
It is known that any counter-example to their question (if one exists)
must fail to be
finitely based in a particularly strange way.  In this paper we show that 
the ``inherently nonfinitely based" algebras constructed by 
Lawrence and Willard from group actions
%to answer two questions from logic
%of A. Wro\'{n}ski and D. Pigozzi 
do \emph{not} fail to be finitely based in this
particularly strange way, and so do not provide a counter-example to the
question of Eilenberg and Sch\"{u}tzenberger.  
As a corollary, we give the first known examples of inherently nonfinitely based ``automatic algebras" constructed from group actions.
\end{abstract}

\maketitle

%%%%%%%%%%%%%%%%%%%%%%%%%%
%% Body

\section{Introduction}
The problem motivating the work presented here concerns finite algebras, as understood in universal algebra, and the identities satisfied in them.  A finite algebra $\m a$
(we will always assume our algebras have just finitely many basic operations)
is said to be \emph{finitely based} (\FB) if the set of identities valid in 
$\m a$ can be axiomatized by some finite subset;
if this is not the case, then $\m a$  is \emph{nonfinitely based} (\NFB).
While many finite algebras of general interest,  including
all finite groups \cite{oates-powell} and finite rings \cites{kruse,lvov:FB}, are \FB,
the situation for general finite algebras is more delicate \cite{tarski}.  Even determining which finite semigroups are \FB\ is an enormous, ongoing problem \cite{volkov}.
In this paper we consider two open questions about finite algebras which are \NFB\ in a particularly
strong way, and solve them in two special classes of algebras 
that arise from finite group actions.

To explain the questions, let us consider a finite algebra $\m a$.
For each $n\geq 1$
we let $\calV(\m a)^{(n)}$ denote the 
class of all models of the set of all $n$-variable identities valid in $\m a$,
so $\calV(\m a)^{(1)}\supseteq \calV(\m a)^{(2)}\supseteq\cdots$, and we let
$\calV(\m a)=\bigcap_{n=1}^\infty \calV(\m a)^{(n)}$.  
Thus $\calV(\m a)$ is the class of all models of the set of all identities
valid in $\m a$; $\calV(\m a)$ is called the \emph{variety generated by} $\m a$.
By a theorem of Birkhoff, $\m a$ is \FB\ if and only if $\calV(\m a)= \calV(\m a)^{(n)}$ for some $n\geq 1$.
Furthermore, $\calV(\m a)$ is
\emph{locally finite}: every finitely generated subalgebra of a 
member of $\calV(\m a)$ is finite.  
So if it happens that $\calV(\m a)^{(n)}$ is \emph{not} locally finite,
for every $n\geq 1$, then $\m a$ must be \NFB.  
We say that $\m a$ is \emph{inherently} \NFB\ (\INFB)
if it satisfies the hypothesis of the previous sentence.
Many known \NFB\ finite algebras are actually \INFB\ 
\cites{murskii,perkins,sapir,shallon,park,jezek,shift,isaev,
law-wil,kea-wil,tarski}.

A notion nominally stronger even than \INFB\ arises if we require finite witnesses
to nonlocal finiteness.
Let $\scrK$ be a class of algebras closed under forming subalgebras and
arbitrary direct products.  If $\scrK$ is not locally finite,
then for some $d\geq 1$ there must exist in $\scrK$ an infinite, 
$d$-generated
algebra.  A stronger property, which may or may not hold, is that for some 
$d\geq 1$ there exist in
$\scrK$ arbitrarily large \emph{finite} $d$-generated algebras. 
Guided by this observation, we say that a finite algebra $\m a$ is \INFB\ 
\emph{with finite witnesses}, or \emph{in the finite sense} (\INFBfin),
if for every $n\geq 1$ there exists $d\geq 1$ such that
$\calV(\m a)^{(n)}$ contains arbitrarily large finite $d$-generated algebras.
%It is not known whether \INFBfin\ is strictly stronger than \INFB:
%
%\medskip\noindent\textsc{\INFBfin\ Problem} \cite{eqcomplexity}:
%
%Does there exist a finite algebra which is \INFB\ but not \INFBfin?

%\medskip
This brings us to the work of Eilenberg and Sch\"{u}tzenberger from 1976.
Motivated by applications to the theory of automata, they considered
finite monoids and asked, for a finite monoid $\m m$, whether the existence
of a finite set of identities whose class of \emph{finite}
models coincides with the class of \emph{finite} members of $\calV(\m m)$
is sufficient to imply that $\m m$ is \FB.  They also noted that their question
can be posed for arbitrary finite algebras, not just monoids:

\medskip\noindent\textsc{Eilenberg-Sch\"{u}tzenberger Problem} \cite{eilenberg}:

If $\m a$ is a finite algebra for which
there exists a finite set of identities whose class of finite models
coincides with the class of finite members of $\calV(\m a)$,
does it follow that $\m a$ is \FB?

\medskip
As Sapir noted in his positive solution to the \Eil\ in the 
case of semigroups \cite{sapir}, if a finite semigroup is a counter-example to the \Eil, then it must be \INFB\ but not \INFBfin.  
His reasoning works generally: any counter-example to the \Eil\ 
must be \INFB\ but not \INFBfin.  This motivates the following
problem, first stated by McNulty et al in 2008:

\medskip\noindent\textsc{\INFBfin\ Problem} \cite{eqcomplexity}:

Does there exist a finite algebra which is \INFB\ but not \INFBfin?

\medskip
Both the \Eil\ and the \INFBfin\ Problem are open.  A reasonable strategy for approaching both problems is to identify 
known classes of \INFB\ algebras and  prove that they are all \INFBfin.
This is how Sapir answered the \Eil\ affirmatively for semigroups, 
using his classification of \INFB\ semigroups \cite{sapir:infb}.
Another source of known \INFB\ algebras is
provided by the so-called \emph{shift automorphism method} \cite{shift}; McNulty
et al \cite{eqcomplexity} proved that every algebra which can be
proved to be INFB by this method is \INFBfin.
A third sporadic source of \INFB\ algebras is the family of
5-dimensional nonassociative bilinear $K$-algebras constructed by Isaev
\cite{isaev}.  These have been studied in a recent preprint by McNulty and Willard
\cite{mcn-wil}
and shown there to be \INFBfin.

In this paper we study another sporadic class of known \INFB\
algebras for which these questions are not already resolved: the algebras
constructed from group actions
by Lawrence and Willard in \cite{law-wil}.
These algebras are constructed as follows.  Let $\alpha$
be a faithful action of a finite group $G$ on a finite set $S$, written $(g,s) \mapsto gs$.
Define an algebra $\Alg GS\alpha^\ast$
with universe $G\times S$ and two basic operations:
\begin{enumerate}
\item
A unary operation $f$ given by $f((g,s)) = (g,gs)$.
\item
A binary operation $d$ given by $d((g_1,s_1),(g_2,s_2)) = (g_1,s_2)$.
\end{enumerate}
In essence, the binary operation allows one to view
$\Alg GS\alpha^\ast$ as a 2-sorted algebra with sorts $G$ and $S$ respectively,
and whose only operation is the action $G\times S\ra S$.  In particular,
although \emph{elements} of $G$ are present and can act on elements
of $S$, the group \emph{operation} of $G$ is not available for use in identities.
Lawrence and Willard \cite{law-wil} proved that (i) if $G$ is nilpotent, then
$\Alg GS\alpha^\ast$ is \FB, while (ii) if $G$ is not nilpotent,
then $\Alg GS\alpha^\ast$ is \INFB.  The main result of our
paper is

\begin{thm} \label{main}
If $G$ is not nilpotent, then $\Alg GS\alpha^\ast$ is \INFBfin.
\end{thm}

Hence no algebra $\Alg GS\alpha^\ast$ arising from a group action in this way
is a counter-example to the \Eil.

Another natural way to model a group action without the group operation 
is as an ``automatic algebra."  Given nonempty sets $G,S$ and a function
$\sigma$ assigning to each $a \in G$ a partial self-map $\sigma_a$ on
$S$ (that is, a function $\sigma_a:U_a\ra S$ for some $U_a\subseteq S$), the \emph{automatic algebra}
determined by $(G,S,\sigma)$ is an algebra $\Auto GS\sigma$
whose universe is the disjoint
union of $G$, $S$ and $\{0\}$ where $0$ plays the role of a default element, and which
has one binary operation $\cdot$ given by
\[
x\cdot y = \left\{\begin{array}{cl}
\sigma_x(y) &\mbox{if $x \in G$ and $y \in \dom(\sigma_x)\subseteq S$}\\
0 & \mbox{otherwise.}
\end{array}\right.
\]
In particular, if $\alpha$ is an action of a group $G$ on a set $S$, we let $\Auto GS\alpha$ denote the automatic algebra $\Auto GS\sigma$ where  each
$\sigma_a$ is the total map $S\ra S$ given by $\sigma_a(s) = as$.

A few small automatic algebras are known to be \INFB\ \cites{kea-wil,boozer,eqcomplexity}, either explicitly or implicitly via
the shift automorphism method.  A limitation of the shift automorphism method is that it doesn't play well with elements (of an algebra) which ``act" as permutations (on other elements of
the algebra).  In particular, the shift automorphism method cannot
be applied to any automatic algebra arising from a group action.
Using our methods presented here, we can overcome this limitation:

\begin{thm} \label{second}
Let $\alpha$ be a faithful action of a finite group $G$
on a finite set $S$.  If $G$ is not nilpotent,
 then $\Auto GS\alpha$ is \INFBfin.
\end{thm}

As in \cite{law-wil}, our proofs ultimately rest on the existence of a family of finitely generated groups with special properties. The innovation here is the fact that the finitely generated groups used in \cite{law-wil}, which are infinite,  have arbitrarily large finite homomorphic images.

Here is an overview of the rest of the paper.
In Section~\ref{improve} we construct the finitely 
presented groups which underpin our theorems, and prove their
needed properties.  In Section~\ref{2sort} we construct and analyze
the 2-sorted versions of the algebras considered in Theorem~\ref{main}.
In Section~\ref{convert} we convert the 2-sorted algebras to their 1-sorted
avatar described above and prove Theorem~\ref{main}.  In Section~\ref{adjoin}
we show how Theorem~\ref{main} can be extended 
to actions of certain semigroups with zero.  Then in Section~\ref{auto} we
prove Theorem~\ref{second}.  The paper concludes with some open problems.

\section{Improving the group construction from \cite{law-wil}} \label{improve}

Recall \cite{neumann}
that $\calA_n\calA_m$ is the variety of all extensions of abelian
groups of exponent dividing $n$ by abelian groups of exponent dividing $m$.
Given $n\geq 2$ and primes $p,q$, let $\calF_{n,p,q}$ denote the 
set of all groups $G$ such that
\begin{enumerate}
\item[$(\ast$)]
There exists $X\subseteq G$ such that $G=\langle X\rangle$, $|X|=n+1$,
and every subgroup of $G$ generated by a proper
subset of $X$ is in $\calA_p\calA_q$.
\end{enumerate}

A key result from \cite{law-wil} is the construction,
for any $n\geq 2$ and any primes $p,q$, of
an infinite group in $\calF_{n,p,q}$.
In this section we improve the construction in order to show that
$\calF_{n,p,q}$ contains arbitrarily large finite groups.
For this we need the following two facts about finite nilpotent groups.

\begin{lm} \label{lm:H1}
If $G$ is a nilpotent group generated by finitely many elements of finite
order, then $G$ is finite.
\end{lm}

\begin{lm} \label{lm:H2}
For each prime $p$ there exist arbitrarily large finite $p$-groups generated
by two elements of order $p$.
\end{lm}

A proof of Lemma~\ref{lm:H1} can be found in \cite[Theorem 3.9(iii)]{suzuki};
see also \cite{jpvee} for an elementary proof.
We are indebted to Eamonn O'Brien, who pointed
us to the following proof of Lemma~\ref{lm:H2}.

\begin{proof}[Proof of Lemma~\ref{lm:H2}]
Let $G=C_p\ast C_p$ be the free product of two cyclic groups of order $p$.
By a theorem of Nielsen \cite{nielsen}, a proof of which can be found in 
\cite{lyndon}, the derived subgroup $G'$ is free of rank $r=(p-1)^2$.  Clearly
$G/G'\cong C_p\times C_p$.  Now fix $n\geq 1$.  There exits a characteristic
subgroup $N$ of $G'$ such that $G'/N \cong C_{p^n}\times\cdots\times C_{p^n}$ 
($r$ factors).  Then $N\lhd G$ and $G/N$ is a finite $p$-group of order 
$p^{nr+2}$.  Finally, $G/N$ inherits from $G$ the property that
it is generated by two elements of order $p$.
\end{proof}

Recall that the \emph{left-normed higher commutators} in a group are defined
by $[x_1,x_2] = x_1^{-1}x_2^{-1}x_1x_2$ and $[x_1,\ldots,x_n] =
[[x_1,\ldots,x_{n-1}],x_n]$ for $n>2$, and that a group $G$ is nilpotent of 
class $\leq c$ if and only if $[x_1,\ldots,x_{c+1}]$ is identically equal to 1 in $G$.
For each prime $p$ and integer $c\geq 1$ let $H_{p,c}$ be group
presented using generators $a_1,a_2$ and the following relations:
\begin{align*}
a_i^p &:i=1,2\\
[a_{i_0},a_{i_1},\ldots,a_{i_c}] &: (i_0,i_1,\ldots,i_c) \in \{1,2\}^{c+1}.
\end{align*}
Then $H_{p,c}$ is nilpotent of class $\leq c$ (see
\cite[Claim 1]{law-wil}), and since $H_{p,c}$ is generated by
two elements of finite order, we get that $H_{p.c}$ is finite by
Lemma~\ref{lm:H1}.  On the other hand, every (finite) nilpotent group
generated by two elements of order $p$ is nilpotent of some class and thus
is a quotient of some $H_{p,c}$;
thus Lemma~\ref{lm:H2} implies $\lim_{c\ra\infty}|H_{p,c}|=\infty$ for each
prime $p$.

\begin{prp}
For each $n\geq 2$ and all primes $p,q$, $\calF_{n,p,q}$
contains arbitrarily large finite groups.
\end{prp}

\begin{proof}
Fix $n,p,q$ as in the statement of the Proposition, and fix $\ell>0$.
We will construct a finite group in $\calF_{n,p,q}$ of order at least $\ell$.

Let $(V,+,0)$ denote a vector space of dimension $n$ over
the $q$-element field, and let $B=\{e_1,\ldots,e_n\}$ be a basis.
For each $i\in \{1,\ldots,n\}$
let $V_i=\Span(B\setminus \{e_i\})$.

Choose $c>0$ large enough so that $|H_{p,c}|\geq \ell$.
Define $G_{V,p,c}$ to be the group
presented using $V$ as the set of generators, and using
the following relations:
\begin{align*}
v^p &: v \in V,\\
[v_0,v_1,\ldots,v_c] &: v_0,\ldots,v_c \in V,\\
[v,w] &: v,w \in V~\mbox{and}~v-w \in V_1 \cup \cdots \cup V_n.
\end{align*}
It follows from \cite[Claim 1]{law-wil} that $G_{V,p,c}$ is nilpotent of class
$\leq c$, and since it is generated by finitely many elements of order $p$,
$G_{V,p,c}$ is finite by Lemma~\ref{lm:H1}.
Because there exist $v,w \in V$ with $v-w\not\in V_1\cup\cdots\cup V_n$,
we have that $H_{p,c}$ is a retract of $G_{V,p,c}$ and hence $|G_{V,p,c}|\geq \ell$.

Now we define an action of $V$ on $G_{V,p,c}$.  For each $x \in V$ define
$\phi_x:V\ra V$ by $\phi_x(v)=v+x$.
Observe that for any $x,v,w \in V$ we have $v-w \in V_1\cup \cdots \cup V_n$
if and only if
$\phi_x(v)-\phi_x(w) \in V_1\cup\cdots\cup V_n$.  
Hence $\phi_x$ extends to an automorphism $\phi_x^\ast$
of $G_{V,p,c}$.  Moreover, the map $\phi^\ast:x\mapsto \phi_x^\ast$ is a group homomorphism
from $V$ to $\Aut G_{V,p,c}$.

Let $G = G_{V,p,c}\rtimes_{\phi^\ast} V$ be the semidirect product 
of $G_{V,p,c}$ by $V$ with respect to $\phi^\ast$.  Clearly $G$ is finite and
$|G|\geq |G_{V,p,c}|\geq \ell$, and it remains to show $G \in \calF_{n,p,q}$.
To distinguish an element $v \in V$ from its image as a generator of $G_{V,p,c}$, 
we shall
denote the latter by $[v]$.  Then the set $X=\{e_1,\ldots,e_n,[0]\}$
generates $G$.
Suppose $Y$ is a subset of $X$ of size $n$.  If $Y=\{e_1,\ldots,e_n\}$
then $\langle Y\rangle = V \in \calA_q\subseteq \calA_p\calA_q$.  Otherwise,
there exists $i$ such that $e_i\not\in Y$.  Then $\langle Y\rangle =H_i \rtimes
V_i$ where $H_i$ is the subgroup of $G_{V,p,c}$ generated by
$\set{[v]}{$v \in V_i$}$.  By design, $H_i$ is abelian and so $H_i \in\calA_p$.
Since
$V_i \in \calA_q$, we have $\langle Y\rangle \in \calA_p\calA_q$ as required.
\end{proof}

\section{Two-sorted algebras} \label{2sort}

Recall that a 
\emph{signature for $2$-sorted algebras} is a function $\tau$ whose domain
is a set of operation symbols and which, for each symbol in its domain,
assigns an 
expression of the form
$i_1\times\cdots\times i_n\ra j$ where $n\geq 0$ and
$i_1,\ldots,i_n,j \in \{1,2\}$.  
This expression is called the \emph{type} of the symbol.  A
\emph{$2$-sorted algebra} in the signature $\tau$ is a structure 
$\m a=(A_1,A_2;\calF)$
where $A_1$ and $A_2$ are sets (the \emph{universes}) and $\calF$ is a set
of finitary operations indexed by the symbols in the domain of
$\tau$, subject to the requirement that if the type of a symbol is
$i_1\times\cdots\times i_n\ra j$ then the corresponding operation must be a
function from
$A_{i_1}\times \cdots\times A_{i_n}$ to $A_j$.  

The standard notions of
subalgebras, products, homomorphisms, congruences,
quotient algebras, terms, and identities are easily extended from
ordinary (1-sorted) algebras to 2-sorted algebras 
(see e.g.\ \cites{bir-lip,goguen,tarlecki,taylor,mck-val}).
The only subtlety arises around the question of whether or not to admit
algebras with one or more empty universe.  In this paper we will only have
need to consider 2-sorted algebras in which both universes are nonempty;
these are called \emph{everywhere nonempty} in \cite{tarlecki}.
If $K$ is a class of everywhere nonempty 2-sorted algebras in the same
signature $\tau$, then we let $\calV(K)$ denote the closure of $K$ under
products, everywhere nonempty subalgebras, and homomorphic images.
If $\m a=(A_1,A_2;\calF)$ is a 2-sorted algebra and $\varnothing \ne
U_i\subseteq A_i$ for $i=1,2$, then we say that $(U_1,U_2)$ \emph{generates}
$\m a$ if the only subalgebra $\m b=(B_1,B_2;\calF)$ of $\m a$ with
$U_i\subseteq B_i$ for $i=1,2$ is $\m a$ itself.  
If $n_1,n_2\geq 1$ then we say 
that $\m a$ is $(n_1,n_2)$-\emph{generated} if $\m a$ is generated by some
$(U_1,U_2)$ with $0<|U_i|\leq n_i$ for $i=1,2$.

For the remainder of this paper, we fix $\tau$ to consist
of one binary operation symbol $s$ of type $1\times 2\ra 2$.  
Let $\Vtau$ denote the class of all everywhere-nonempty
algebras in this signature.

\begin{df}
Suppose $G$ is a group, $S$ is a set, and $\alpha:G\times S\ra S$ is a
faithful left action of $G$ on $S$.  We define $\Alg GS\alpha$ to be the
two-sorted algebra $(G,S;\alpha)$ in $\Vtau$.
%with universes $G$ and $S$ and whose only 
%operation is $\alpha$.
\end{df}

In \cite[\S4]{law-wil} the authors worked with the 2-sorted algebra $(S,G;\alpha)$
where the type of $\alpha$ is now $2\times 1\ra 1$; the choice
of how to order the universes of $\Alg GS\alpha$ is a matter of taste and
makes no material difference to the results in this paper.

We need the following definition and
lemma guaranteeing that certain algebras are in $\calV(\Alg GS\alpha)$.

\begin{df}
If $H$ is a group and $r\geq 1$, then we let $H^{\otimes r}$ denote the
disjoint union of $r$ copies of $H$, and let $\m l(H,r)$ denote
$\Alg H{H^{\otimes r}}\lambda$ where $\lambda$ is the action of $H$ on $H^{\otimes r}$ by
left multiplication.
\end{df}

\begin{lm} \label{lm:otimes}
Suppose $G$ is a group, $\alpha:G\times S\ra S$ is a faithful left action,
and $H \in \calV(G)$.
%Let $H^{\oplus r}$ be the
%disjoint union of $r$ copies of $H$, and let
%$\beta$ be the action of $H$ on $H^{\oplus r}$ by left multiplication.
Then $\m l(H,r) \in \calV(\Alg GS\alpha)$ for all $r\geq 1$.
\end{lm}

\begin{proof}
For each $g \in G$ define $\Delta(g)$ to be the constant map in $G^S$
with value $g$, and define $R(g)$ to be the map in $S^S$ given by
$R(g)(s)=\alpha(g,s)$.  It can be easily checked that $(\Delta,R)$ is
an embedding of $\m l(G,1)$ into $\Alg GS\alpha^S$, so $\m l(G,1)
\in \calV(\Alg GS\alpha)$.

Next note that for any groups $H,K$ we have:
\begin{enumerate}
\item
$\m l(H^X,1) = \m l(H,1)^X$ for any set $X$.
\item
If $K\leq H$ then $\m l(K,1) \leq \m l(H,1)$.
\item
If $N\lhd K$ and $\theta_K$ is the corresponding congruence of $K$,
then $(\theta_N,\theta_N)\in \Con(\m l(K,1))$ and
$\m l(K,1)/(\theta_N,\theta_N)\cong \m l(K/N,1)$.
\end{enumerate}
It follows from these facts that if $H \in \calV(G)$ then $\m l(H,1) \in \calV(
\m l(G,1))$ and hence $\m l(H,1) \in \calV(\Alg GS\alpha)$.

Finally, given $r\geq 1$, note that $H \in \calV(G)$ implies $H^r \in \calV(G)$
so $\m l(H^r,1) \in \calV(\Alg GS\alpha)$.  
One can show that $\m l(H,r)$
embeds into $\m l(H^r,1)$. (Hint: $H^r$ has a subgroup
isomorphic to $H$; identify $H^{\otimes r}$ with the union of $r$ distinct
right cosets of this subgroup.)
So $\m l(H,r) \in \calV(\Alg GS\alpha)$ as well.
\end{proof}

\begin{thm}\label{thm:B}
Suppose $G$ is a finite group, $\alpha:G\times S\ra S$ is a faithful left
action, and $G$ is not nilpotent.  Then for every $n\geq 2$ and $\ell>0$
there exists a 2-sorted algebra $\m b(n,\ell)
\in \Vtau$ 
%in the same signature as $\m a(\alpha)$,
satisfying:
\begin{enumerate}
\item \label{B:it1}
$\m b(n,\ell)$ is $(n+1,1)$-generated.
\item \label{B:it2}
Both universes of $\m b(n,\ell)$ are finite and nonempty.
\item \label{B:it3}
The second universe of $\m b(n,\ell)$ has size $\geq \ell$.
\item \label{B:it4}
Each everywhere-nonempty $(n,n)$-generated subalgebra of $\m b(n,\ell)$ belongs to 
$\calV(\Alg GS\alpha)$.
\end{enumerate}
\end{thm}

\begin{proof}
By \cite[Theorem 2.4]{law-wil}, we can choose and fix primes $p,q$ with
$p\ne q$ and such that $\calA_p\calA_q\subseteq \calV(G)$.
Given $n\geq 2$ and $\ell>0$, choose $P \in \calF_{n,p,q}$ so that $P$ is
finite and $|P|\geq \ell$.  Choose $X=\{x_0,x_1,\ldots,x_n\}\subseteq P$
witnessing the fact that $P\in \calF_{n,p,q}$; that is, $P=\langle X\rangle$
and every subgroup
of $P$ generated by a proper subset of $X$ is in $\calA_p\calA_q$.
Recall that $\m l(P,1)$ is the algebra $(P,P;\lambda)$ where $\lambda$ is 
left multiplication.
Now let $\m b(n,\ell)$ be the subalgebra of $\m l(P,1)$
with universes $X$ and $P$ respectively.

Clearly $\m b(n,\ell)$ is generated by $(X,\{1\})$, proving \eqref{B:it1}.
Items \eqref{B:it2} and \eqref{B:it3} are obviously true.  
So it remains to prove \eqref{B:it4}.
Suppose $\m c=(C_1,C_2;\lambda|_{C_1\times C_2})$ is an $(n,n)$-generated subalgebra
of $\m b(n,\ell)$ with $C_1,C_2\ne\varnothing$.  
Note that $C_1$ is a proper subset
of $X$.  Let $H$ be the subgroup of $P$ generated by $C_1$.
By hypothesis, $H \in \calA_p\calA_q$ so $H \in \calV(G)$.
$C_2$ is closed under the action of $H$ by left multiplication, so $C_2$
is the union of some number, $r$, of left cosets of $H$ in $P$.  Thus
$\m c$ embeds into $\m l(H,r)$.
As $\m l(H,r) \in \calV(\Alg GS\alpha)$ by Lemma~\ref{lm:otimes}, we get $\m c \in \calV(\Alg GS\alpha)$.
\end{proof}

\section{Conversion to 1-sorted algebras}  \label{convert}

In this section we exploit a general categorical 
conversion of the variety of all everywhere-nonempty $k$-sorted
algebras in a fixed signature to a variety of 
%(nonempty) 
1-sorted algebras, due ultimately to Barr \cite[Theorem 5]{barr}
(see \cite[\S5]{goguen} for a fuller account, 
and \cite[chapter 11]{mck-val} for explicit details).
The resulting 1-sorted variety is 
defined only up to term equivalence, meaning that the conversion
is really between 
multi-sorted and 1-sorted clones (for universal algebraists) or algebraic
theories (for category theorists).  To avoid having to explain the machinery of
clones or algebraic theories, we will keep this presentation concrete by 
explaining the conversion in the special case of the 
signature considered in the previous section, 
choosing one specific presentation of the clone of the 1-sorted variety which
is produced. 

So again let 
$\tau$ denote the 2-sorted signature consisting of exactly one binary 
operation $s$ of sort $1\times 2\ra 2$, and recall
that $\Vtau$ is the class of all
everywhere-nonempty 2-sorted algebras in the signature $\tau$.
Let $\tau^\ast$ be the signature of 1-sorted algebras consisting of a
binary operation $d$ and a unary operation $f$.  Let $\Vtau^\ast$ denote
the variety of 1-sorted algebras in the signature $\tau^\ast$ 
axiomatized by the following identities:
\begin{align*}
d(x,x) &\approx x\\
d(d(x,y),d(z,w)) &\approx d(x,w)\\
d(f(x),y) &\approx d(x,y).
\end{align*}
For every $\m a=(A_1,A_2;s) \in \Vtau$, define a 1-sorted algebra
$\m a^\ast$ in the signature $\tau^\ast$ as follows: the universe is
$A_1\times A_2$, and the operations $d,f$ are given by
\begin{align*}
d((a_1,a_2),(b_1,b_2)) &= (a_1,b_2)\\
f((a_1,a_2)) &= (a_1,s(a_1,a_2)).
\end{align*}

\begin{prp}[essentially Barr \cite{barr}] \label{prp:barr}
$\m a\mapsto \m a^\ast$ is (the object map of) a category equivalence from $\Vtau$ to $\Vtau^\ast$.
\end{prp}

Indeed, given an algebra $\m c \in \Vtau^\ast$, define binary relations
$E_1,E_2$ on $C$ by $aE_1b$ if and only if
there exists $c \in C$ with $d(a,c)=d(b,c)$,
and $aE_2b$ if and only if there exists $c \in C$ with $d(c,a)=d(c,b)$.  One can use
the defining identities of $\Vtau^\ast$ to show that $E_1$ and $E_2$ 
are equivalence relations
and the map $C\ra C/E_1 \times C/E_2$ given by $c \mapsto (c/E_1,c/E_2)$ is 
a bijection with inverse given by $(a/E_1,b/E_2) \mapsto d(a,b)$.
In particular, if $aE_1a'$ and $bE_2b'$ then $d(a,b)=d(a',b')$.
Now define the 2-sorted algebra $\m c^\square \in \Vtau$ to have universes
$C/E_1$ and $C/E_2$ and operation $s:C/E_1\times C/E_2\ra C/E_2$ given by
$s(a/E_1,b/E_2) = f(d(a,b))/E_2$.  One can verify $(\m a^\ast)^\square
\cong \m a$
naturally for all $\m a \in \Vtau$, and $(\m c^\square)^\ast \m c$
naturally for all $\m c \in \Vtau^\ast$.

Additionally, the following facts are easily verified.
\begin{enumerate}
\item[(F1)] \label{F1}
For all $\m a \in \Vtau$ and $n>0$, $\m a$ is $(n,n)$-generated if and only if
$\m a^\ast$ is $n$-generated.
\item[(F2)]
For all $\m a \in \Vtau$, if $\m c\leq \m a^\ast$ then $\m c=\m b^\ast$
for some $\m b\leq \m a$.
\item[(F3)]
For all $\m a,\m b \in \Vtau$, $\m b \in \calV(\m a)$ 
if and only if
$\m b^\ast \in \calV(\m a^\ast)$.
\end{enumerate}

Now we can prove the main result of this paper.

\medskip\noindent\textbf{Theorem~\ref{main}.}
\emph{Suppose $G$ is a finite group, $\alpha:G\times S\ra S$ is a faithful
left action, and $G$ is not nilpotent.  Then $\Alg GS\alpha^\ast$ is
\INFBfin.
}

\begin{proof}
Let $V=\calV(\Alg GS\alpha^\ast)$.
Fix $n\geq 2$.
We will show that $V^{(n)}$ contains arbitrarily
large finite $d$-generated algebras where $d=n+1$.
Fix $\ell>0$.  
Let $\m b =\m b(n,\ell)=(B_1,B_2;s)$ be the 
2-sorted algebra from Theorem~\ref{thm:B}.  Then $\m b \in \Vtau$,
both $B_1$ and $B_2$ are finite, $|B_1|\geq 1$, and $|B_2|\geq \ell$.
Consider the 1-sorted image $\m b^\ast \in \Vtau^\ast$ of $\m b$.
We have $\m b^\ast$ is finite and $|B^\ast|=|B_1|\cdot|B_2|\geq \ell$.
Moreover, (F1) implies $\m b^\ast$ is $d$-generated because $\m b$ is $(d,1)$-generated.

It remains to prove $\m b^\ast \in V^{(n)}$.
Because $V^{(n)}$ is defined by the $n$-variable identities valid in $V$, it suffices to show that every $n$-generated subalgebra of $\m b^\ast$ is in $V$
(see \cite[Lemma~7.15]{ALVINv2}).
Suppose
$\m c$ is an $n$-generated subalgebra of $\m b^\ast$.  Then
$\m c=\m d^\ast$ for some $(n,n)$-generated subalgebra $\m d$ of $\m b$,
by (F2) and (F1).  By Theorem~\ref{thm:B}, $\m d \in \calV(\Alg GS\alpha)$, 
so $\m c \in \calV(\Alg GS\alpha^\ast)=V$ by (F3).  This proves $\m b^\ast \in V^{(n)}$.
\end{proof}

%%%%%%%%%%%%%%%%%%%%%%%%%%%%%%%%%%%%%%%%%%
\section{Adjoining zero} \label{adjoin}

The construction of $\Alg GS\alpha^\ast$ does not require a group action; any action of
one finite set on another
%$\alpha:G\times S\ra S$ of a finite semigroup $G$ on a finite set 
will do.  
%One could ask for which finite  faithful semigroup actions $\alpha$ is $\m a(\alpha)^\ast$ \INFB\ (or \INFBfin).
In this section we extend our main theorem to actions by finite semigroups consisting of a nonnilpotent group with a zero element adjoined.
More precisely, let $G$ be a group and let $\alpha$ be an action of
$G$ on a set $S$.  Let $G^0$ denote the semigroup obtained from $G$ by adding a zero
element 0, let $S^0$ denote the set obtained from $S$ by adding a 
new element 0, and let $\alpha^0$ denote the extension of $\alpha$ to a semigroup action of $G^0$ on $S^0$ satisfying
$\alpha^0(0,y) = \alpha^0(x,0)=0$ for all $x \in G^0$
and all $y \in S^0$.  We will show that if $G$ and $S$ 
are finite, $G$ is nonnilpotent, and $\alpha$ is faithful,
then $\Alg{G^0}{S^0}{\alpha^0}$ is \INFBfin.

Given any 2-sorted algebra $\m b=(B_1,B_2;s)
\in \Vtau$, let $\m b^0$ denote the 2-sorted algebra
$(B_1^0,B_2^0;s^0)$, where $B_i^0$ is the disjoint union of $B_i$ with
$\{0\}$ and $s^0$ is the extension of $s$ given by 
\begin{align*}
{s}^0(b,b') &= \left\{\begin{array}{cl}
{s}(b,b') & \mbox{if $b \in B_1 $ and $b' \in B_2$}\\
0 & \mbox{otherwise.}
\end{array}\right.
\end{align*}
Note in particular that if $G,S,\alpha$ are as before, then
$\Alg{G^0}{S^0}{\alpha^0} = \Alg GS\alpha^0$.  
For readability, we will also denote
$\Alg GS\alpha^0$ by $\Algzero GS\alpha$.

The following facts are easily verified.
\begin{enumerate}
\item[($F'_1$)]
For $\m a, \m b \in \Vtau$, if $\m b\leq \m a$ then $\m b^0\leq \m a^0 $.
\item[($F'_2$)]
For $\m a_i \in \Vtau$, $ (\prod_{i} \m a_i )^0 $ embeds into $ \prod_{i} \m a_i^0  $.
\item[($F'_3$)]
For $\m a,\m b \in \Vtau$, any homomorphism ${h}:\m a\ra \m b$ has an extension $h^0:\m a^0\ra \m b^0$ such that ${h^0}(\m a^0)=(h(\m a))^0$. Furthermore, if $\m c\leq \m a$, then ${h^0}(\m c^0)\cong(h(\m c))^0$. 
\item[($F'_4$)] 
For all $\m a \in \Vtau$ and $\m b\leq \m a^0$, $\m b$ is $(X_1,X_2)$-generated if and only if $\m b^- = (B_1^-,B_2^-;s)$ is $(X_1^-,X_2^-)$-generated, where $B_i^- = B_i\setminus \{0\}$ and $X_i^- = X_i\setminus \{0\}$. 
\end{enumerate}

\begin{lm} \label{L1}
Fix $\m a \in \Vtau$.
If $\m b \in \calV(\m a)$, then $\m b^0 \in\calV(\m a^0)$.
\end{lm}
\begin{proof}
By $(F_1'), (F_2')$ and $(F_3')$, the class $\set{\m b \in \Vtau}{$\m b^0 \in \calV(\m a^0)$}$ contains $\m a$ and is closed under $\Ho$, $\Su$ and $\Pd$.
\end{proof}

\begin{thm} \label{thm:zero}
Let $G$ be a finite nonnilpotent group and let $\alpha$ be a faithful action
of $G$ on a finite set $S$.  
%Let $G^0$ denote the semigroup obtained by
%adding a zero element to $G$.  Also let $S^0$ be the disjoint union of $S$ with $0$, and let $\alpha^0$ be the action of $G^0$ on $S^0$ given by
%\begin{align*}
%\alpha^0(g,s) &= \left\{\begin{array}{cl}
%\alpha(g,s) & \mbox{if $g \in G$ and $s \in S$}\\
%0 & \mbox{otherwise.}
%\end{array}\right.
%\end{align*}
Then $\Algzero GS\alpha^\ast$ is \INFBfin.
\end{thm}

\begin{proof}
%Observe that $\Alg {G^0}{S^0}{\alpha^0}=\Alg GS\alpha^0$.
Let 
$V=\calV(\Alg GS\alpha)$, $V_0 = \calV(\Algzero GS\alpha)$, 
and $V_0^\ast = \set{\m a^\ast}{$\m a \in V_0$}$,
and 
observe that $\calV(\Algzero GS\alpha^\ast)$ is the closure of $V_0^\ast$ under isomorphisms by Proposition~\ref{prp:barr} 
and Fact (F3) from the previous section.
Fix $n\geq 2$.
We will show that $(V_{0}^{\ast})^{(n)}$ contains arbitrarily
large finite $d$-generated algebras where $d=n+2$.
Now for fixed $n\geq 2$ and $\ell>0$  Theorem~\ref{thm:B} gives an algebra $\m b(n,\ell)\in \Vtau$ satisfying four properties.  Let $\m b=\m b(n,\ell)$ and consider
$\m b^0$. Clearly
\begin{enumerate}
\item[$(1)'$]
$\m b^0$ is $(n+2,2)$-generated.
\item[$(2)'$]
Both universes of $\m b^0$ are finite and nonempty.
\item[$(3)'$]
The second universe of $\m b^0$ has size $>\ell$.
\end{enumerate}
Moreover we can show:
\begin{enumerate}
\item[$(4)'$]
Each everywhere-nonempty $(n,n)$-generated subalgebra of $\m b^0$ is in $V_0$.
\end{enumerate}
Indeed, 
suppose $\m d$ is an everywhere-nonempty $(n,n)$-generated subalgebra of $\m b^0$.
If $D_1=\{0\}$ or $D_2=\{0\}$ then $\m d \in V_0$ is easily
verified.  So assume that $D_1\ne\{0\}$ and $D_2\ne \{0\}$.
%It remains to prove that $\m d \in V^{\alert{0}}$.
%
%\medskip\noindent\textsc{Case 1}: $D_1=\{0\}$ or $D_2=\{0\}$. Trivial.
%
%\medskip\noindent\textsc{Case 2}: Otherwise.  
Let $D_i^- = D_i\setminus \{0\}$.
Then 
%$D_i^-\ne\varnothing$ for $i=1,2$. Note that 
$(D_1^-,D_2^-)$ are universes of an 
everywhere-nonempty $(n,n)$-generated
 subalgebra $\m d^-$ of $\m b$. Thus $\m d^- \in V$ 
 by one of the properties of $\m b$ and 
 so $(\m d^-)^0 \in V_0$ by Lemma \ref{L1}. 
 Further, $\m d\leq (\m d^-)^0$, which implies $\m d \in V_0$, proving $(4)'$.

%so $\m c \in \calV((\m a(\alpha^0)^\ast)=V^0$ by (F3).  

Consider the 1-sorted image $(\m b^0)^\ast \in \Vtau^\ast$ of $\m b^0$.
We have $(\m b^0)^\ast$ is finite, $|(B^0)^\ast|\geq \ell$ and $(\m b^0)^\ast$ is $n+2$-generated. 
Suppose $\m c$ is an $n$-generated subalgebra of $(\m b^0)^\ast$. % Then
By fact (F1) from the previous section,
$\m c=\m d^\ast$ for some $(n,n)$-generated subalgebra $\m d$ of $\m b^0$.
Then
$\m d \in V_0$ by $(4)'$, so $\m c \in V_{0}^{\ast}$.
%by fact (F3).
This proves that
$(\m b^0)^\ast \in (V_{0}^{\ast})^{(n)}$.
\end{proof}

%\medskip
%Note that in particular if $G=S_3$ and $\mathsf{\lambda}$ \alert{is the} faithful representation of $S_3$ as the set of permutations on $\{1,2,3\}$, then the algebras $\m a(\alert{\lambda}) $ and $\m a(\alert{\lambda}^0)$ described in Theorem\ref{main} and Theorem \ref{thm:zero} have $18$ and $28$ elements respectively and further both the algebras are inherently nonfinitely based in the finite sense.   

%%%%%%%%%%%%%%%%%%%%%%%%%%%%
\section{Application to automatic algebras} \label{auto}

%The semigroup actions $\alpha^0:G^0\times S^0\ra S^0$ considered in the previous section are examples of  \emph{semigroup-with-zero actions}.  
%These are actions
 %$\alpha:G\times S\ra S$ where the semigroup $G$ is required to have
 %a zero element $0$ satisfying $0x=x0=x$ for all $x \in G$,
 %the set $S$ is \emph{pointed}, meaning that it has a privileged element also denoted $0$,
 %and the action satisfies $\alpha(x,0)=\alpha(0,y)=0$ for all $x \in G$
 %and $y \in S$.  (What is \emph{not} required is that $x \in G\setminus \{0\}$ and $y\in S\setminus \{0\}$ imply $\alpha(x,y)\ne 0$.)

%Given a semigroup-with-zero action $\alpha$ of a semigroup-with-zero $G$ on a pointed set $S$, 
%there is another natural way to model $\alpha$:
%as an \emph{automatic
%algebra}.  This is the 1-sorted algebra 
%$\aut(\alpha)$ which is obtained by identifying the
%zero elements of $G$ and $S$, assuming that $G$ and $S$ are
%otherwise disjoint (i.e., $G\cap S=\{0\}$), and taking
%$G\cup S$ to be the universe of $\aut(\alpha)$;
%5 and then defining one binary operation $\cdot$ on the universe by
%\[
%x\cdot y = \left\{\begin{array}{cl}
%\alpha(x,y) & \mbox{if $x \in G$ and $y \in S$}\\
%0 & \mbox{otherwise.}
%\end{array}\right.
%\]
Recall that an automatic algebra is determined by 
%
%the definition of ``automatic algebra" in the
%introduction: given 
nonempty sets $G,S$ and a function
$\sigma$ assigning to each $a \in G$ a partial self-map $\sigma_a$ on
$S$, and is denoted by $\Auto GS\sigma$.
%%
%
%, the \emph{automatic algebra}
%determined by $(G,S,\sigma)$ is the algebra $\Auto GS\sigma$
%whose universe is the disjoint
%nion of $G$, $S$ and $\{0\}$, and which
%has one binary operation $\cdot$ given by
%\[
%x\cdot y = \left\{\begin{array}{cl}
%\sigma_x(y) &\mbox{if $x \in G$ and $y \in \dom(\sigma_x)\subseteq S$}\\
%0 & \mbox{otherwise.}
%\end{array}\right.
%\]
We also use the notation $\Auto GS\alpha$ 
% for the automatic algebra 
when $\alpha$ is an action of a group $G$ on a set $S$.
%\alert{If $\alpha$ is an action of a group $G$ on a set $S$, then
%$\Auto GS\alpha$ denotes $\Auto GS\sigma$ where
%$\sigma_a(s)=as$ for all $a \in G$ and $s \in S$.}

In this section we show how the proof of Theorem~\ref{thm:zero} can be modified to prove the second main result of this paper.

\medskip\noindent\textbf{Theorem~\ref{second}.}
\emph{Suppose $G$ is a finite nonnilpotent group and $\alpha$ is a
faithful action of $G$ on a finite set $S$.  Then
$\Auto GS\alpha$ is \INFBfin.
}

\medskip
Before starting the proof, we introduce some notation and one
fact which will be useful.    Let $X$ be an infinite
set of variables.  Given a word $w$ over $X$ 
and a variable $y$ we define
a term $[w]y$ in the signature of automatic algebras as follows: 
\begin{itemize}
\item
If $\varepsilon$ is the empty word then $[\varepsilon]y = y$.
\item
If $w=xw'$ with $x\in X$, then $[w]y = x\cdot [w']y$.
\end{itemize}
For example, $[x_1x_2x_3]y = x_1\cdot (x_2\cdot (x_3\cdot y))$.

Now consider the 2-sorted signature
$\tau$ with its single operation $s$ of type $1\times 2\ra 2$.
If $w$ is a word over $X$ and $y$ is a variable not in $X$, then
we can reinterpret $[w]y$ as a term of sort 2 in the signature $\tau$ by:
\begin{itemize}
\item
reinterpreting the variables in $X$ as having sort 1, and the 
variable $y$ as having sort 2;
\item
replacing each occurrence of $\cdot$ by $s$.
\end{itemize}
Denote this reinterpretation by $[w]y^\sharp$.
For example, $[x_1x_2x_3]y^\sharp = s(x_1,s(x_2,s(x_3,y)))$
where now $x_1,x_2,x_3$ are variables of sort 1 and $y$ is a variable of sort 2.

Let $\Auto GS\sigma$ be an arbitrary automatic algebra. 
We use $\mathbf{0}$ as an abbreviation for the term $z\cdot z$.  Note that $\mathbf{0}$ is identically equal to 0 in $\Auto GS\sigma$.  
Let $\Alg {G^0}{S^0}\sigma$ denote the 2-sorted algebra $(G^0,S^0;s)$ in $\Vtau$ whose operation $s$ is the restriction of the operation $\cdot$ of $\Auto GS\sigma$ to $G^0\times S^0$.

Now define the following sets of identities
in the signature of automatic algebras:
\begin{align*}
\Delta &= \{\mathbf{0}\cdot x \approx \mathbf{0},~ x\cdot\mathbf{0} \approx \mathbf{0},~x\cdot x\approx \mathbf{0},
~(x\cdot y)\cdot z \approx \mathbf{0}\}\\
&\qquad{}\cup~ \set{[x_0x_1\cdots x_n]x_0 \approx \mathbf{0}}{$n\geq 0$},\\
\Psi(G,S,\sigma) &= \set{[w]y\approx [w']y}{$[w]y^\sharp \approx [w']y^\sharp$ is identically true in $\Alg {G^0}{S^0}\sigma$},\\
\Psi_0(G,S,\sigma) &= \set{[w]y \approx \mathbf{0}}{$[w]y^\sharp$  is identically equal to 0 in $\Alg {G^0}{S^0}\sigma$}.
\end{align*}

We need the following fact, which is a consequence of the
``Core Theorem" in John Boozer's PhD thesis.

\begin{fact}[Boozer \cite{boozer}] \label{fact:boozer}
Let $\Auto GS\sigma$ be an automatic algebra.
Then 
\[
\Delta \cup \Psi(G,S,\sigma) \cup \Psi_0(G,S,\sigma)
\]
is a basis for $\calV(\Auto GS\sigma)$.  
%That is,
%\begin{enumerate}
%\item
%Every identity in $\Delta \cup \Psi\cup \Psi_0$ is identically
%true in $\aut(\alpha)$, and
%\item
%Every identity which is identically true in $\aut(\alpha)$ is
%a logical consequence of $\Delta\cup\Psi\cup \Psi_0$.
%\end{enumerate}
\end{fact}

We also need one more well-known fact, due to Birkhoff.

\begin{fact}[Birkhoff \cite{birkhoff}] \label{fact:birk}
If $\m a$ is a finite algebra in a finite signature and
$n>0$, then $\calV(\m a)^{(n)}$ is finitely based.
\end{fact}

\begin{proof}[Proof of Theorem~\ref{second}]
Recall that $G$ is a finite nonnilpotent group, $\alpha$ is a faithful action of $G$ on the finite set $S$, and we wish to show
that the automatic algebra $\Auto GS\alpha$ is \INFBfin.  
%\alert{Also recall that $\Auto GS\alpha=\Auto GS\sigma$
%where $\sigma_a(s)=\alpha(a,s)$ for all $a \in G$ and $s \in S$.}

Fix $n\geq 2$.  By 
Fact~\ref{fact:birk}, there exists a finite basis $\Sigma_n$ for
$\calV(\Auto GS\alpha)^{(n)}$.  
By Fact~\ref{fact:boozer} and 
the compactness theorem, there exists a finite subset
$\Sigma_n' \subseteq 
\Psi(G,S,\sigma) \cup \Psi_0(G,S,\sigma)$ such that every identity
in $\Sigma_n$ is a logical consequence of $\Delta\cup\Sigma_n'$.  In fact, because 
$\alpha(x,y)\ne 0$ for all $x \in G$ and $y \in S$, it follows
that $\Psi_0(G,S,\sigma)=\varnothing$, so we have $\Sigma_n'\subseteq \Psi(G,S,\sigma)$.
Let $k$ be the maximum number of variables occurring in an identity in $\Sigma_n'$.  Let $d=k+4$.  We will show that $\calV(\Auto GS\alpha)^{(n)}$ has arbitrarily large finite
$d$-generated members.

%Let $V_0 = \calV(\m a(\alpha^0))$.

Recall the 2-sorted algebra $\Algzero GS\alpha$ defined in the previous section, and let $V_0 = \calV(\Algzero GS\alpha)$.  Note that the 2-sorted algebra $\Alg{G^0}{S^0}\sigma$ defined before Fact~\ref{fact:boozer} is in our current context just $\Algzero GS\alpha$.

Fix $\ell>0$.  
Recall that in
the proof of Theorem~\ref{thm:zero} we constructed a 2-sorted algebra
$\m b^0=\m b(k,\ell)^0$ satisfying these properties:
\begin{enumerate}
\item[$(1)'$]
$\m b^0$ is $(k+2,2)$-generated.
\item[$(2)'$]
Both universes of $\m b^0$ are finite and nonempty.
\item[$(3)'$]
The second universe of $\m b^0$ has size $>\ell$.
\item[$(4)'$]
Each everywhere-nonempty $(k,k)$-generated subalgebra of $\m b^0$ is in $V_0$.
% $\calV(\Alg GS\alpha^0)$.
\item[$(5)'$]
Let the universes of $\m b^0$ be $(B_1,B_2)$, and let its
operation be $s$.  Then each $B_i$ contains an element denoted
$0$ such that $s(0,y)=s(x,0)=0$ for all $x \in B_1$ and all
$y \in B_2$.
\end{enumerate}

%Clearly $|B_1|\leq k+2$.
Define an algebra $\m c$ in the signature of automatic
algebras as follows: arrange that $B_1\cap B_2=\{0\}$, let the universe of
$\m c$ be $B_1\cup B_2$, and define the operation $\cdot$ on
$B_1 \cup B_2$ by
\[
x\cdot y = \left\{\begin{array}{cl}
s(x,y) & \mbox{if $x \in B_1$ and $y \in B_2$}\\
0 & \mbox{otherwise.}
\end{array}\right.
\]
Note that if $(X_1,X_2)$ generates $\m b^0$ then
$X_1\cup X_2$ generates $\m c$, so $\m c$ is
$k+4$-generated.  Clearly $\m c$ is finite and
$|C|\geq |B_2|>\ell$.  It remains to show
that $\m c$ is in $\calV(\Auto GS\alpha)^{(n)}$,
or equivalently, that $\m c\models \Sigma_n$.  It
suffices to show $\m c \models \Delta\cup \Sigma_n'$.
Clearly $\m c$ satisfies every identity in $\Delta$, so
we need only check the identities in $\Sigma_n'$.  Suppose
$[w]y\approx [w']y$ is an identity in $\Sigma_n'$.  
By construction, $[w]y \approx [w']y$ is in $\Psi(G,S,\sigma)$,
which means $y$ does not occur in $w$ or $w'$ and $\Alg {G^0}{S^0}\sigma\models [w]y^\sharp \approx [w']y^\sharp$, or equivalently,
$\Algzero GS\alpha\models [w]y^\sharp \approx [w']y^\sharp$.
Observe that this last fact implies that $w$ and $w'$ contain exactly the same variables; otherwise if a variable $x$
were to occur in $w$, say, but not in $w'$, then
we could falsify $[w]y^\sharp\approx [w']y^\sharp$ in $\Algzero GS\alpha$ by an
assignment sending $x\mapsto 0$ and all other variables to
nonzero elements of the appropriate sorts.

We now prove that $\m c \models [w]y\approx [w']y$.
Let $x_1,\ldots,x_r$ be a list of the distinct variables that
occur in $w$ (equivalently in $w'$).  
Note that $r<k$ by our choice of $k$.
Let $x_i\mapsto a_i$, $y\mapsto b$ be
an assignment of values in $C$ to the variables in
$\{x_1,\ldots,x_n\}\cup \{y\}$.  Assume that this assignment
falsifies the identity $[w]y\approx [w']y$.  This can only
happen if
$a_1,\ldots,a_r \in B_1\setminus \{0\}$, $b \in B_2
\setminus \{0\}$, and the same assignment falsifies $[w]y^\sharp \approx [w']y^\sharp$ in $\m b^0$.  This
falsification is then witnessed in the subalgebra of
$\m b^0$ generated by $(\{a_1,\ldots,a_r\},\{b\})$.
But every $(k,k)$-generated subalgebra of $\m b^0$
is in $V_0$ by $(4)'$, 
so cannot falsify
$[w]y^\sharp \approx [w']y^\sharp$ as $V_0\models [w]y^\sharp \approx [w']y^\sharp$.
\end{proof}

\section{Summary and questions}

Given a finite nonnilpotent group acting faithfully on a finite set, we have described three finite algebras that
capture the group action: $\Alg GS\alpha^\ast$, $\Algzero GS\alpha^\ast$,
and $\Auto GS\alpha$.  The first of these was invented in
\cite{law-wil} and shown there to be inherently nonfinitely based (\INFB).  The second is a simple variation of the first.  The
third is an example of an automatic algebra to which the 
shift automorphism method does not apply, and was not previously
known to be \INFB. In this paper we 
showed that all three algebras are inherently nonfinitely based in the finite sense (\INFBfin).  Thus none of these algebras
can be a counterexample to the \Eil.

In particular, if $G=S_3$ and $\mathsf{\alpha}$ is the faithful representation of $S_3$ as the set of permutations on $\{1,2,3\}$, then the algebras $\Alg GS\alpha^\ast$,  $\Algzero GS\alpha^\ast$ and $\Auto GS\alpha$ have $18$, $28$, and 10 elements respectively, and all are inherently nonfinitely based in the finite sense.  
%Each of these algebras is proved here for the first time to be inherently nonfinitely based in the finite sense. 

We end by posing three problems.

\begin{prb}
For which finite semigroups $G$ does there exist a finite $G$-action $\alpha$ such that 
$\Alg GS\alpha^\ast$ is inherently nonfinitely based (\INFB)?
\end{prb}

\begin{prb}
Does the \Eil\ have a positive answer for algebras of
the form $\Alg GS\alpha^\ast$ where $\alpha$ is an action
of a finite semigroup $G$ on a finite set $S$?
\end{prb}

\begin{prb}
Does the \Eil\ have a positive answer for automatic algebras?
\end{prb}

\bibliography{infqbfinbib}

\end{document}